\documentclass{article}

\usepackage{a4wide}

\begin{document}

\title{Optimization of Dengue Epidemics: a test case\\ 
with different discretization schemes\thanks{This is a preprint
of a paper whose final and definite form has been published in 
\emph{American Institute of Physics Conf. Proc. {\bf 1168} (2009), no.~1, 1385--1388. 
{\tt DOI:10.1063/1.3241345}}}}

\author{Helena Sofia Rodrigues\\
\texttt{sofiarodrigues@esce.ipvc.pt}\\
School of Business Studies\\
Viana do Castelo Polytechnic Institute, Portugal
\and
M. Teresa T. Monteiro\\
\texttt{tm@dps.uminho.pt}\\
Department of Production and Systems\\
University of Minho, Portugal
\and
Delfim F. M. Torres\\
\texttt{delfim@ua.pt}\\
Department of Mathematics\\
University of Aveiro, Portugal}

\date{}

\maketitle


\begin{abstract}
The incidence of Dengue epidemiologic disease has grown in recent decades. 
In this paper an application of optimal control in Dengue epidemics is presented. 
The mathematical model includes the dynamic of Dengue mosquito, the affected persons, 
the people's motivation to combat the mosquito and the inherent social cost of the disease, 
such as cost with ill individuals, educations and sanitary campaigns. 
The dynamic model presents a set of nonlinear ordinary differential equations. 
The problem was discretized through Euler and Runge Kutta schemes, and solved 
using nonlinear optimization packages. The computational results 
as well as the main conclusions are shown.

\bigskip


\smallskip


\smallskip

\end{abstract}

\maketitle


\section{Introduction}

Dengue is a mosquito mostly found in tropical and sub-tropical climates worldwide, 
mostly in urban and semi-urban areas. It can provoque a mosquito-borne infection, 
that causes a severe flu-like illness, and sometimes a potentially lethal complication 
called dengue haemorrhagic fever, and about 40\% of the world's population are now at risk.

The aim of this paper is to present an attempt to apply quantitative methods in the optimization of
investments in the control of Epidemiologic diseases, in order to
obtain a maximum of benefits from a fixed amount of financial
resources. This model includes the dynamic of the growing of the mosquito, 
but also the efforts of the public management to motivate the population to break the reproduction 
cycle of the mosquitoes by avoiding the accumulation of still water in open-air 
recipients and spraying potential zones of reproduction.

The paper is organized as follows. Next section presents the dynamic model for dengue epidemics, 
where the variables, parameters and the control system are defined. Then, the numerical implementation 
and the strategies used to solve the problem are reported. 
Finally, the numerical results are presented and some conclusions are taken.


\section{Dynamic model}

The model described in this paper is based on the model proposed in \cite{Caetano2001}.

The notation used in the mathematical model is as follows.

\medskip

\begin{tabular}{lp{10cm}}
\multicolumn{2}{l}{\textbf{State Variables:}}\\
$x_{1}(t)$ & density of mosquitoes; \\
$x_{2}(t)$ & density of mosquitoes carrying the virus; \\
$x_{3}(t)$ & number of persons with the disease; \\
$x_{4}(t)$ & level of popular motivation to combat mosquitoes (goodwill). \\
\end{tabular}

\medskip

\begin{tabular}{lp{10cm}}
\multicolumn{2}{l}{\textbf{Control Variables:}}\\
$u_{1}(t)$ & investments in insecticides; \\
$u_{2}(t)$ & investments in educational campaigns. \\
\end{tabular}

\medskip

\begin{tabular}{lp{12.5cm}}
\multicolumn{2}{l}{\textbf{Parameters:}}\\
$\alpha_{R}$ & average reproduction rate of mosquitoes; \\
$\alpha_{M}$ & mortality rate of mosquitoes; \\
$\beta$ & probability of contact between non-carrier mosquitoes and affected persons; \\
$\eta$ & rate of treatment of affected persons; \\
$\mu$ & amplitude of seasonal oscillation in the reproduction rate of mosquitoes; \\
$\rho$ & probability of persons becoming infected; \\
$\theta$ & fear factor, reflecting the increase in the population's willingness 
to take actions to combat the mosquitoes as a consequence of the high prevalence 
of the disease in the specific social environment; \\
$\tau$ & forgetting rate for goodwill of the target population; \\
$\varphi$ & phase angle to adjust the peak season for mosquitoes; \\
$\omega$ & angular frequency of the mosquitoes's proliferation cycle, corresponding to 52 weeks period;\\
$P$ & population in the risk area (usually normalized to yield $P=1$); \\
$\gamma_D$ & the instantaneous costs due to the existence of affected persons;\\
$\gamma_F$ & the costs of each operation of spraying insecticides;\\
$\gamma_E$ & the cost associated to the instructive campaigns. \\
\end{tabular}

\bigskip

The Dengue epidemic can be modeled by the following nonlinear time-varying state equations. 
Equation (\ref{x1}) represents the variation of the density of mosquitoes per unit time 
to the natural cycle of reproduction and mortality ($\alpha_R$ and $\alpha_M$), 
due to seasonal effects $\mu sin(\omega t+\varphi)$ and to human interference $- x_4(t)$ and $u_1(t)$:

\begin{equation}\label{x1}
\frac{dx_{1}}{dt}=\left[\alpha_{R}\left(1-\mu sin(\omega t+\varphi)\right)-\alpha_M - x_4(t)\right]x_1(t)-u_1(t).
\end{equation}

Equation (\ref{x2}) expresses the variation of the density $x_2$ of mosquitoes carrying the virus. 
The term $\left[\alpha_{R}\left(1-\mu sin(\omega t+\varphi)\right)-\alpha_M - x_4(t)\right]x_2(t)$ 
represents the rate of the infected mosquitoes and $\beta \left[ x_1(t)-x_2(t)\right]x_3(t)$ 
represents the increase rate of the infected mosquitoes due to the possible contact between 
the non infected mosquitoes $x_1(t)-x_2(t)$ and the number persons with disease denoted by $x_3$:

\begin{equation}\label{x2}
\frac{dx_{2}}{dt}=\left[\alpha_{R}\left(1-\mu sin(\omega t+\varphi)\right)-\alpha_M 
- x_4(t)\right]x_2(t)+\beta \left[ x_1(t)-x_2(t)\right]x_3(t)-u_1(t).
\end{equation}

The dynamics of the infectious transmission is presented in equation (\ref{x3}). 
The term $-\eta x_3(t)$ represents the rate of cure and $\rho x_2(t)\left[P-x_3(t)\right]$ 
represents the rate at which new cases spring up. The factor $\left[P-x_3(t)\right]$ 
is the number of persons in the area, that are not infected:

\begin{equation}\label{x3}
\frac{dx_{3}}{dt}= -\eta x_3(t)+\rho x_2(t)\left[P-x_3(t)\right].
\end{equation}

Equation (\ref{x4}) is a model for the level of popular motivation (or goodwill) 
to combat the reproductive cycle of mosquitoes. Along the time, the level of people's 
motivation changes and, as consequence, it is necessary to invest in educational 
campaigns designed to increase consciousness of the population under risk by 
a proper understanding of the determinants involved with specific disease. 
The expression $-\tau x_4(t)$ represents the decay of the people's motivation with time, 
due to forgetting. The expression $\theta x_3(t)$ represents the natural sensibilities 
of the public due to increase in the prevalence of the disease.

\begin{equation}\label{x4}
\frac{dx_{4}}{dt}= -\tau x_4(t)+\theta x_3(t)+u_2(t).
\end{equation}

The goal of the problem is to minimize the cost functional
\begin{equation}\label{cost}
J\left[u_1(.),u_2(.)\right]=\int_{0}^{t_f}\{\gamma_D x_{3}^{2}(t)
+\gamma_F u_{1}^{2}(t)+\gamma_E u_{2}^{2}(t)\}dt.
\end{equation}

This functional includes the social costs related to the existence of ill persons, 
$\gamma_D x_{3}^{2}(t)$, the recourses needed for spraying of insecticides operations, 
$\gamma_F u_{1}^{2}(t)$, and for educational campaigns, $\gamma_E u_{2}^{2}(t)$. 
The model for the social cost is based in the concept of goodwill 
explored by Nerlove and Arrow \cite{Nerlove1962}.

Due to computational issues, the optimal control problem (\ref{x1})-(\ref{cost}) 
that is in the Lagrange form, was converted into an equivalent Mayer form.
Hence, using a standard procedure to rewrite the cost functional \cite{Lewis1995}, 
the state vector was augmented by an extra component $x_5$,

\begin{equation}\label{newx5}
\frac{dx_5}{dt}=\gamma_D x_{3}^{2}(t)+\gamma_F u_{1}^{2}(t)+\gamma_E u_{2}^{2}(t)
\end{equation}

\noindent leading to the equivalent terminal cost problem of minimizing

\begin{center}
\begin{tabular}{l}
$I[x_5(.)]=x_5(t_f)$
\end{tabular}
\end{center}

\noindent with given $t_f$, subject to the control 
system (\ref{x1})-(\ref{x4}) and (\ref{newx5}).


\section{Numerical implementation}

The simulations were carried out using the following normalized numerical values:
$\alpha_{R}=0.20$, $\alpha_{M}=0.18$, $\beta=0.3$, $\eta=0.15$, $\mu=0.1$, $\rho=0.1$, 
$\theta=0.05$, $\tau=0.1$, $\varphi=0$, $\omega=2\pi/52$, $P=1.0$, $\gamma_D=1.0$, 
$\gamma_F=0.4$, $\gamma_E=0.8$, $x_1(0)=1.0$, $x_2(0)=0.12$, $x_3(0)=0.004$, and $x_4(0)=0.05$. 
These values are available in the paper \cite{Caetano2001}. The final time used was $t_f=52$ weeks.

To solve this problem it was necessary to discretize the problem. Two methods were selected: 
a first order, the Euler's scheme, and a Runge Kutta's sheme of second order \cite{Betts}. 
In both cases, it is assumed that the time $t=nh$ moves ahead in uniform steps of length $h$. 
If a differential equation is written like $\displaystyle\frac{dx}{dt}=f(t,x)$, it is possible 
to make a convenient approximation of this. In the Euler's scheme the update is given by
\begin{center}
\begin{tabular}{l}
$x_{n+1}\simeq x_n+hf(t_n,x_n),$\\
\end{tabular}
\end{center}

\noindent while in the Runge Kutta's method is

\begin{center}
\begin{tabular}{l}
$x_{n+1}\simeq x_n+\frac{h}{2}\left[f(t_n,x_n)+f(t_{n+1},x_{n+1})\right].$\\
\end{tabular}
\end{center}

This approximation $x_{n+1}$ of $x(t)$ at the point $t_{n+1}$ has an error depending 
on $h^2$ and $h^3$, for the Euler and Runge Kutta methods, respectively. 
This discretization process transforms the dengue epidemics problem into a standard 
nonlinear optimization problem (NLP), with an objective function and a set of nonlinear constraints. 
This NLP problem was codified, for both discretization schemes, in the AMPL modelling language \cite{AMPL}.

Two nonlinear solvers with distinct features were selected to solve the NLP problem: the Knitro and the Snopt. 
The first one \cite{Knitro} is a software package for solving large scale mathematical optimization 
problems based mainly on the Interior Point (IP) method. Snopt \cite{Snopt} 
uses the SQP (Sequential Quadratic Programming) philosophy, 
with an augmented Lagrangian approach combining a
trust region approach adapted to handle the bound constraints. 
The NEOS Server \cite{NEOS} platform was used as interface with
both solvers.


\section{Computational results}

Table~\ref{resultados} reports the results for both solvers, for each discretization method 
using three different discretization steps ($h=0.5, 0.25,  0.125$), rising twelve numerical 
experiences. The columns \# var. and \# const. mean de number of variables and constraints, 
respectively. The next columns refer to the performance measures -- number of iterations 
and total CPU time in seconds (time for solving the problem, for evaluate the objective 
and the constraints functions and for input/output). The computational experiences 
were made in the NEOS server platform - in this way the selected machine 
to run the program remain unknown as well as its technical specifications.

\begin{table}
{\footnotesize
\begin{tabular}{ccc}
\hline
& Euler's method & Runge Kutta's method\\
\hline
\hline
Knitro & \begin{tabular}{|c|c|c|c|c|}
           \hline
           h & \# var. & \# const. & \# iter. &  time (sec.)\\ \hline
           0.5 & 727 & 519 & 113 & 2.090 \\
           0.25 & 1455 & 1039 & 68 & 2.210 \\
           0.125 &  2911 & 2079 & 85 & 7.240 \\
           \hline
         \end{tabular}
& \begin{tabular}{|c|c|c|c|c|}
           \hline
           h & \# var. & \# const. & \# iter. &  time (sec.)\\ \hline
           0.5 & 728 & 520 & 64 & 1.980 \\
           0.25 & 1456 & 1040 & 82 & 5.550 \\
           0.125 & 2912 & 2080 & 70 & 9.740 \\
           \hline
         \end{tabular}
\\ \hline

Snopt & \begin{tabular}{|c|c|c|c|c|}
           \hline
           h & \# var. & \# const. & \# iter. & time (sec.)\\ \hline
           0.5 & 727 & 519 & 175 & 4.07 \\
           0.25 &  1455 & 1039 & 253 & 19.2 \\
           0.125 & 2911 & 2079 & 252 & 105.4 \\
           \hline
         \end{tabular}
& \begin{tabular}{|c|c|c|c|c|}
           \hline
           h & \# var. & \# const. & \# iter. &  time (sec.)\\ \hline
           0.5 & 728 & 520 &  223 & 10.52 \\
           0.25 & 1456 & 1040 & 219 &  39.7 \\
           0.125 & 2912 & 2080  & 420 &  406.67\\
           \hline
         \end{tabular}
\\ \hline
\end{tabular}}
\caption{Numerical results}
\label{resultados}
\end{table}

The optimal value reached was $\approx 3E-03$ for all tests. 
Comparing the general behaviour of the solvers one can conclude that the IP based method (Knitro) 
presents much better performance than the SQP method (Snopt) in terms of the measures used. 
Regarding the Knitro results, one realize that the Euler's discretization scheme has better 
times for $h=0.25$ and $h=0.125$ and similar time for $h=0.5$, when compared to Runge-Kutta's method. 
Another obvious finding, for both solvers, is that the CPU time increases as far as the problem dimension 
increases (number of variables and constraints). With respect to the number of iterations, 
Snopt presents more iterations as the problem dimension increases. However this conclusion 
cannot be taken for Knitro -- in fact, doesn't exist a relation between the problem dimension 
and the number of iterations.  The best version tested was Knitro using Runge-Kutta with $h=0.5$ 
(best CPU time and fewer iterations), and the second one was Knitro with  Euler's method using $h=0.25$. 
An important evidence of this numerical experience is that it is not worth the reduction 
of the discretization step size because no significative advantages are obtained.


\section{Conclusions}

We solved successfully an optimal control problem by direct methods using nonlinear 
optimization software based on IP and SQP approaches. The effort of the implementation 
of higher order discretization methods brings no advantages. The reduction 
of the discretization step and consequently  the increase of the number of variables 
and constraints doesn't improve the performance with respect to the CPU time 
and to the number of iterations. We can point out the robustness 
of both solvers in spite of the dimension problem increase.


\section*{Acknowledgments}

The first author is grateful to
the Portuguese Foundation for Science and Technology 
(FCT) for the PhD Grant SFRH/BD/33384/2008.




\begin{thebibliography}{9}

\bibitem{Betts}
J.~Betts, \emph{Practical Methods for Optimal Control Using Nonlinear Programming},
SIAM: Advances in Design and Control, 2001.

\bibitem{Knitro}
R.H.~Byrd, and J.~Nocedal, and R.A.~Waltzy,
``Knitro: An Integrated Package for Nonlinear Optimization,''in
\emph{Large-Scale Nonlinear Optimization}, G. di Pillo and M. Roma, eds, Springer-Verlag, 2006, pp. 35-59.

\bibitem{Caetano2001}
M.~Caetano, and T.~Yoneyama, Optimal and sub-optimal control in Dengue epidemics,
\emph{Optimal Control Applications and Methods} \textbf{22}, 2001, 63--73.

\bibitem{AMPL}
R.~Fourer, and D.M.~Gay, and B.W.~Kernighan,
\emph{AMPL: A Modeling Language for Mathematical Programming},
Duxbury Press / Brooks/Cole Publishing Company, 2002.

\bibitem{Snopt}
P.E.~Gill, and W.~Murray, and M.A.~Saunders,
SNOPT: An SQP Algorithm for Large-Scale Constrained Optimization,
\emph{SIAM}, \textbf{47}-1, pp 99--131, 2005.

\bibitem{Lewis1995}
F.L.~Lewis, and V.L.~Syrmos,
\emph{Optimal Control (2nd ed)},
Wiley: New  York, 1995.

\bibitem{NEOS}
NEOS: \emph{www-neos.mcs.ang.gov/neos, July 2009}

\bibitem{Nerlove1962}
M.~Nerlove, and K.J.~Arrow,
Optimal advertising policy under dynamic conditions,
\emph{Economica} \textbf{42}-114, 129--142, 1962.

\end{thebibliography}
\end{document}